%
%
\input amstex.tex
\input amsppt.sty
\documentstyle{amsppt}
\input pictex
\NoBlackBoxes
\magnification=\magstep1
\advance\vsize-0.5cm\voffset=-0.5cm\advance\hsize1cm\hoffset0cm
\font\svtnrm cmr17 scaled \magstep2

\topmatter
\title Applications of controlled surgery in dimension 4: Examples \endtitle
\author Friedrich Hegenbarth and Du\v san Repov\v s \endauthor
\address
Department of Mathematics,
University of Milano, Via C. Saldini 50, Milano, Italy 02130.
e-mail: friedrich.hegenbarth\@mat.unimi.it \endaddress
\address
Institute for Mathematics, Physics and Mechanics,
University of Ljubljana, Jadranska 19, Ljubljana, Slovenia 1001.
e-mail: dusan.repovs\@fmf.uni-lj.si \endaddress

\subjclass Primary 57RXX, Secondary 55PXX
\endsubjclass
\keywords 
Generalized manifold, ANR, Poincar\'e duality, 
$\varepsilon$-$\delta$--surgery, controlled, Poincar\'e complex,
Quinn index, disk theorem, knot group, spine
\endkeywords

\abstract
The validity of Freedman's disk theorem is known to
depend only on the fundamental group. It was conjectured that it
fails for nonabelian free fundamental groups.
If this were true
then surgery theory 
would work 
in dimension four.
Recently, Krushkal 
and Lee proved a surprising result 
that surgery theory
works for a large special class of 4-manifolds
with free nonabelian fundamental groups.
The goal of this paper is to show that
this also holds for other fundamental groups
which are not known to be good, and that it is best
understood using controlled surgery theory
of 
Pedersen--Quinn--Ranicki.
We consider some examples of 4-manifolds
which have the fundamental group
either
of a closed aspherical surface 
or of a 
3-dimensional knot space.
A more general theorem is stated in the appendix.
\endabstract
\endtopmatter

\document

\head {$\S$ 1. Introduction} \endhead

The purpose of this paper is to study
$4$--dimensional surgery problems by means of
controlled surgery.
The usual higher dimensional
surgery procedure breaks down in dimension four
since framed $2$--spheres can generically  only be immersed
in a $4$--manifold (whereas for surgery on them
one would require
embeddings). To get an
embedding one uses the {\sl Whitney trick}. 
Its basic ingredient
is the existence of {\sl Whitney disks} along which
pairs of intersection points with
opposite algebraic intersection number can be
cancelled. If one finds these Whitney disks,
surgery can be completed
provided that the
Wall obstruction vanishes. The celebrated
Disk Theorem of Freedman asserts that (see \cite{Fr}):

\item{(1)}
The existence of Whitney disks in a
$4$--manifold  $M^4$ depends only on the fundamental
group of $M^4$. If they exist then $\pi_1 (M^4)$ is called a {\sl good}
fundamental group.
\item{(2)}
The (large) class of good fundamental groups
includes the trivial group and $\Bbb Z$.
(see also \cite{Fr-Qu}, \cite{Fr-Tei}, \cite{Kru-Qu}).

\noindent
It has been conjectured that nonabelian
free groups are not good.
Nevertheless, the following surprising
result was proved by Krushkal and Lee (\cite{Kru-Lee}):

\proclaim{Theorem 1.1} 
Let $X$ be a $4$--dimensional Poincar\' e complex with a
free nonabelian fundamental group, and assume that
the intersection form on $X$ is extended from
the integers. Let $f : M \to X$ be a degree one normal map, 
where $M$ is a closed $4$--manifold.
Then the vanishing of the Wall obstruction implies
that $f$ is normally bordant to
a homotopy equivalence $f' : M' \to X$.
\endproclaim

Whenever the intersection form is extended from the integers, it follows that
$X$ is homotopy equivalent to
the connected sum $P \# M'$, where $M'$ is simply
connected and $P$ can be assumed to be
a finite sum
$\mathop\#\limits _1\limits^r S^1 \times S^3$
(see \cite{He-Re-Sp}).
Recall that there exist
Poincar\' e $4$--complexes with free fundamental group and
intersection form which is not extended from $\Bbb Z$
(see \cite{Ham-Tei}, \cite{He-Pic}).

By means of
this example we see that
surgery can be completed for a
large class of $4$--manifolds (or Poincar\' e complexes)
with fundamental group $\pi$, which
is not good.
This paper will confirm this fact for
other fundamental groups. For instance we shall prove:

\proclaim{Theorem 1.2}
Let $X$ be a spin Poincar\' e $4$--complex and suppose that
it 
has the fundamental group of a 
closed oriented aspherical surface
and that
the intersection form is extended from $\Bbb Z$. Then
any degree one normal map $f : M \to X$ with
vanishing Wall obstruction is normally bordant
to a homotopy equivalence.
\endproclaim

We shall also give other examples, e.g.
4-manifolds having the
fundamental group isomorphic to some special knot group.
Moreover, we shall also recover Theorem 1.1.

The reason why this can happen is that
one can divide
the {\sl global} surgery
problem into smaller pieces
for which the {\sl local}
fundamental groups are
good, i.e. $\{1\}$ or $\Bbb Z$. 
One gets several
local surgery obstructions
which assemble
to give the global surgery obstruction. This
subdivision has to be done in such a way
that the global surgery obstruction already
determines the local ones. More precisely,
the assembly map should be injective.

The subdivision is made according
to a control map $p : X \to B$, where
$B$ is a finite-dimensional compact metric
ANR.
The map $p$ must satisfy the following three
conditions:

\item{(i)}
$p$ is a $UV^1$--map;

\item{(ii)}
if $X$ is a Poincar\' e complex then
$X$ must be a $\delta$--controlled Poincar\' e
complex over $B$, where $\delta > 0$ is smaller
than some $\varepsilon_0 > 0$ which depends only on $B$; and

\item{(iii)}
the assembly map $A : H_4 (B, \Bbb L) \to L_4 (\pi_1 (B))$
is injective.

\noindent
Definitions and more explanations will be given
in Section 2.
Once we have such a control map we
can apply controlled surgery theory to obtain
our results.
Note that the extreme cases, i.e. when

\item{(a)}
either $p = \operatorname{Id} : X \to X = B$;

\item{(b)}
or $p = \operatorname{const} : X \to \{*\} = B$,

\noindent
generically do not work since the case (a) does not
satisfy condition (iii) whereas the case (b) does not satisfy condition (i) above.
It is also obvious that $p : X \to B$
depends not only on $\pi_1(X)$ but also on the
topology of $X$, so one gets solution
of the surgery problem in individual
cases.

\bigskip
{\it Remark.} The condition (i) can be weakened
to:

\item{(i')} $UV^1(\delta)$, for every $0 < \delta < \varepsilon_0$.

\noindent
Another example of a 4-manifold with a knot group fundamental 
group is stated in Theorem 3.2 below.

\vfill\eject


\head {$\S$ 2. Controlled surgery theory} \endhead

To reach our goal we shall 
need to use
the $\varepsilon$--$\delta$ surgery
sequence and to
compare it with the
non-controlled one.
Let $\Bbb L$ denote the $4$--periodic simply
connected surgery spectrum (see \cite{Qu1} and \cite{Ni} for
geometric and \cite{Ra} for algebraic definitions).
For a space $B$ we then have $\Bbb L$--homology
(resp. cohomology) groups denoted by $H_p(B, \Bbb L)$
(resp. $H^p(B, \Bbb L)$). There is a well defined
assembly map $A : H_p(B, \Bbb L) \to L_p(\pi_1(B))$, where
$L_p(\pi_1(B))$ denotes the Wall group of obstructions 
to simple homotopy equivalences.

We shall only
consider the oriented situation. 
Let $X$ be an
$n$--dimensional simple Poincar\' e complex. We
suppose that $X$ admits a degree one normal map
$f_0 : M_0 \to X$, so fixing this we have
an identification of all degree one normal maps into $X$,
modulo a normal cobordism,
with the homotopy set $[X, G/_{TOP}]$ (see \cite{Wa}).

There is a well defined map
$$
\Theta : [X, G/_{TOP}] \to H_n (X, \Bbb L),
$$
which associates to a given degree one normal map
into $X$ its local surgery problems
according to a small dissection of $X$
(see \cite{Ra}). Composition with the assembly
map
$$
\sigma : A \circ \Theta : [X, G/_{TOP}] \to L_n(\pi_1(X))
$$
yields the classical surgery obstruction
map of Wall.

The topological structure set of $X$ consists
of simple homotopy equivalences $f : M^n \to X$,
where $M^n$ is a closed manifold. Two such
homotopy equivalences $f : M \to X$, $f' : M' \to X$
are equivalent if there is a homeomorphism
$g : M \to M'$ such that the diagram


\centerline{
\beginpicture 
\setcoordinatesystem units <0.800mm,0.800mm>
\setplotarea x from -20 to 20, y from -15 to 15
\setsolid
\setplotsymbol ({\fiverm.})
\arrow <2mm> [.2,.4] from 10.000 5.000 to 5.000 -5.000
\arrow <2mm> [.2,.4] from -10.000 5.000 to -5.000 -5.000
\arrow <2mm> [.2,.4] from -10.000 10.000 to 10.000 10.000
\put {$f'$} at 10.000 0.000
\put {$f$} at -10.000 0.000
\put {$g$} at 0.000 13.000
\put {$X$} at 0.000 -10.000
\put {$M'$} at 15.000 10.000
\put {$M$} at -15.000 10.000
\endpicture
}

\noindent
homotopy commutes.
The set of equivalence classes will be
denoted by $\Cal S (X)$ and will be called the 
{\sl topological
structure}
set of $X$.

Any homotopy equivalence determines a degree one normal map
inducing
$
\Cal S (X) \to [X, G/_{TOP}].
$
For $n \ge 5$ and $X$ a simple Poincar\' e complex,
there is an exact ordinary surgery sequence
$$
\Cal S (X) \to [X, G/_{TOP}] \to L_n(\pi_1(X)).
$$

This sequence can be extended to the left by
the Wall realization of obstructions
$$
L_{n+1}(\pi_1(X)) \to \Cal S (X).
$$
For this one has to fix a simple homotopy equivalence
$f_0 : M_0 \to X$, where $M_0$ is a topological manifold and
$\operatorname{dim} X \ge 6$. In the controlled concept there
is a realization of elements in $H_{n+1} (B, \Bbb L)$ giving
a four--term exact sequence which also holds
for $n = 4$. In this paper we do not consider this
part of the sequence so we will not give more
details.

Before we state the $\varepsilon$--$\delta$ surgery sequence we
need some
more definitions.
Let $p : X \to B$ be a control map, $B$ a
(finite-dimensional) compact metric ANR.
Then $p$ is a $UV^1 (\delta)$--map, $\delta > 0$, if
every commutative diagram

\centerline {
\beginpicture 
\setcoordinatesystem units <0.800mm,0.800mm>
\setplotarea x from -20 to 20, y from -20 to 20
\setsolid
\setplotsymbol ({\fiverm.})
\arrow <2mm> [.2,.4] from -14.000 4.000 to -14.000 -5.000
\circulararc 180.000 degrees from -14.000 4.000 center at -15.000 4.000
\arrow <2mm> [.2,.4] from -10.000 10.000 to 10.000 10.000
\arrow <2mm> [.2,.4] from -10.000 -10.000 to 10.000 -10.000
\arrow <2mm> [.2,.4] from 15.000 5.000 to 15.000 -5.000
\put {$K_0$} at -15.000 10.000
\put {$X$} at 15.000 10.000
\put {$K$} at -15.000 -10.000
\put {$B$} at 15.000 -10.000
\put {$\alpha_0$} at 0.000 13.000
\put {$\alpha$} at 0.000 -7.000
\put {$p$} at 18.000 0.000
\endpicture
}

\noindent
where $K$ is a $2$--complex and $K_0 \subset K$ is a subcomplex,
can be completed by a map $\overline{\alpha} : K \to X$
such that $\overline{\alpha} \Big\vert_{K_0} = \alpha_0$ 
and
$d(p \circ  \;{\overline{\alpha}} (u), \alpha (u)) < \delta$ 
for all $u \in K$.
The map $p$ is called a $UV^1$ map if it is a 
$UV^1 (\delta)$ map
for every
$\delta > 0$.
Here,
$d : B \times B \to \Bbb R_+$ denotes the metric on $B$.

Suppose now that $X$ is an $n$--dimensional
Poincar\' e complex, and suppose that
$X$ has a 
simplicial structure. By the Borsuk theorem
the simplicial structure hypothesis can
often be obtained replacing $X$ by a homotopy
equivalent space and then working with it.

A space $X$ is said to be an (oriented) $\delta$--{\sl Poincar\' e complex}
with respect to $p : X \to B$ if

\item{(i)}
for every simplex $\Delta$ of $X$ the diameter
of $p(\Delta) \subset B$ is less than $\delta$; and

\item{(ii)}
there exists a fundamental cocycle $\xi \in C_n (X)$
such that the cap product
$$
\cap \xi : C^k (X) \to C_{n-k} (X)
$$
is a $\delta$--chain equivalence.

The second condition requires geometric module
structure on the $\Lambda$--chain complex
$\{C_k (X) = H_n (\mathop{} \limits \tilde X^{(k)}, \mathop{} \limits \tilde X^{(k-1)})
\;\;\; | \;\;\; k = 0,1,2, \dots \}$,
where $\Lambda = \Bbb Z [\pi_1(X)]$.
We shall not give any more details but will refer
to the literature (see \cite{Ra-Ya1}, \cite{Ra-Ya2}). A more
geometric definition was given in \cite{Qu2}.

If $X$ is a manifold then one obtains by barycentric
subdivision a $\delta$--Poincar\' e structure
for every $\delta > 0$ with respect to $p = \operatorname{Id} : X \to X$,
hence with respect to every $p : X \to B$.

Suppose that $f, g : Y \to X$ are given maps. Then
$f$ is said to be $\delta$--{\sl homotopic} to $g$ if there is a
homotopy $h : Y \times I \to X$ between $f$
and $g$ such that for any $y \in Y$ the
diameter of $\{ph (y, t) \;\;\; | \;\;\; t \in I\} \subset B$
is less than $\delta$.

Moreover, $f : Y \to X$ is called a $\delta$--{\sl homotopy equivalence}
if there exists $g : X \to Y$ and homotopies
$h : X \times I \to X$, $h' : Y \times I \to Y$ between
$f \circ g$ and $\operatorname{Id}_X$ (resp. $g \circ f$ and $\operatorname{Id}_Y$) such that
the diameters of $\{ph (x, t) \;\;\; | \;\;\; t \in I\}$ and
$\{pfh' (y, t) \;\;\; | \;\;\; t \in I\}$ are less than $\delta$ for
all $x \in X$ (resp. $y \in Y$).

\proclaim{Lemma 2.1}
Suppose that
$f : Y \to X$ is a $\delta$--homotopy
equivalence and $Y$ a $\delta'$--Poincar\' e complex
with respect to $p \circ f$. Then $X$ is a
$(\delta' + 2 \delta)$--Poincar\' e complex over $B$.
\endproclaim

This is a useful observation. It follows 
easily from \cite{Ra-Ya1; Proposition 2.3}.
The following is the main theorem of \cite{Pe-Qu-Ra}:

\proclaim{Theorem 2.2}
Let $B$ be a finite--dimensional compact ANR and $X$ a closed
topological $n$--manifold, where $n \ge 4$.
Then there exists an $\varepsilon_0 > 0$,
depending on $B$, such that for every 
$0 < \varepsilon < \varepsilon_0$
there exists $\delta > 0$ such that the following holds: 
If there is a map $p : X \to B$ satisfying the 
$UV^1 (\delta)$ property, then we get
the following controlled surgery exact sequence
$$
H_{n+1} (B, \Bbb L) \to \Cal S_{\varepsilon, \delta} (X, p) \to
[X, G/_{TOP}] \to H_n (B, \Bbb L).
$$
\endproclaim

In fact, the following supplement holds
(however, we shall not need it):

\proclaim{Supplement}
If the map, $p : X \to B$ is only assumed to be a sufficiently small controlled
Poincar\' e complex over $B$ 
(instead of assuming that $X$ is a manifold),
then the surgery sequence
$$
\Cal S_{\varepsilon, \delta} (X,p) \to
[X, G/_{TOP}] \to H_n (B, \Bbb L)
$$
is still exact.
\endproclaim

The controlled structure set $\Cal S_{\varepsilon, \delta} (X,p)$ is
defined as follows: Its elements are represented
by $\delta$--homotopy equivalence $f : M \to X$ 
over $B$, where $M$ is a closed topological $n$--manifold.
Another  $\delta$--homotopy equivalence $g : N \to X$
is said to be
$\varepsilon$--{\sl related} to $(M, f)$ if there exists
a homeomorphism $h : M \to N$, so that


\centerline{
\beginpicture 
\setcoordinatesystem units <0.800mm,0.800mm>
\setplotarea x from -20 to 20, y from -15 to 15
\setsolid
\setplotsymbol ({\fiverm.})
\arrow <2mm> [.2,.4] from 10.000 5.000 to 5.000 -5.000
\arrow <2mm> [.2,.4] from -10.000 5.000 to -5.000 -5.000
\arrow <2mm> [.2,.4] from -10.000 10.000 to 10.000 10.000
\put {$g$} at 10.000 0.000
\put {$f$} at -10.000 0.000
\put {$h$} at 0.000 13.000
\put {$X$} at 0.000 -10.000
\put {$N$} at 15.000 10.000
\put {$M$} at -15.000 10.000
\endpicture
}

\noindent
is $\varepsilon$--homotopy commutative.
This relation is reflexive and symmetric, so it
gives rise to an equivalence relation.
(However, the proof actually shows the
transitivity of the relations above.)
This theorem has a relative version
which we shall not use.

As explained above, we now have
the following commutative diagram for $n = 4$
(assuming the hypotheses of the theorem),

$$
\CD
\Cal S_{\varepsilon, \delta}(X,p) @>>> [X, G/_{TOP}] @>>> H_4 (B, \Bbb L) \\
@VVV	             		     @|		        @VV A V   \\
\Cal S(X) 			@>>> [X, G/_{TOP}] @>>> L_4 (\pi_1(B)). \\
\endCD
$$

Hence in order to solve the $4$--dimensional surgery
problem with target $X$ one needs
a control map $p : X \to B$ satisfying
conditions (i), (ii) and (iii) stated in $\S 1$.
Examples will be given in the
following section.


\head {$\S$ 3. Examples} \endhead

There is another characterization of $UV^1$--maps
which is useful in the applications (see \cite{Dav} or
\cite{BFMW}). A subset $A$ of a space $X$ is said to be
$UV^1$ if for each neighborhood $U$ of $A$ in $X$ there
is another neighborhood $V$ of $A$ with $V \subset U$,
such that the induced map
$
\pi_1 (V) \to \pi_1 (U)
$
is zero for any base point in $V$, and any
two points in $V$ can be connected in $U$.
The following is a special case of the Approximate
lifting theorem (\cite{Dav; p. 126}):

\proclaim{Theorem 3.1} {\rm (\cite{Dav})}
Suppose $X$ is a metric space and $G$ is an upper semicontinuous 
$UV^1$--decomposition of $X$ (i.e.
each member $A \in G$ is
a $UV^1$ subset). Let $B = X/_G$
and
$p : X \to B$.
Then $p$ is a 
$UV^1$
map, i.e. $p$ is a $UV^1(\delta)$ map for every $\delta > 0$.
\endproclaim

As the first example we consider
$X$ as a Poincar\' e $4$--complex with free
nonabelian fundamental group and
$\Bbb Z$--extended $\Lambda$--intersection form, $\Lambda = \Bbb Z [\pi_1(X)]$.
By results from \cite{He-Re-Sp}, $X$ is homotopy
equivalent to $(\mathop\#\limits _1\limits^r S^1 \times S^3) \# M' = M$.
A homotopy equivalence $M \to X$
induces an "isomorphism" of the ordinary short
exact surgery sequences, i.e.
we transform a surgery problem with target $X$
to a surgery problem with target $M$.

\proclaim{Lemma 3.2}
Let $M \# M'$ be the connected sum
of two topological manifolds and $p : M \# M' \to M$ the
map which collapses $M'$ to a point.
If $M'$ is simply connected then $p$ is a $UV^1$ map.
\endproclaim

\demo{Proof}
The map $p$ is the composition of the maps 
$p_2$ and $p_1$. First,
$p_1 : M \# M' \to M \vee M'$ 
is the map  which
collapses
the
$3$--sphere $\sum^3 \subset M \# M'$ to
the base point of the wedge
$M \vee M'$. More precisely,
a bicollar $[-1, 1] \times \sum^3$ is radially
smashed to $D^4 \vee {D'}^4 \subset M \vee M'$, fixing
$\{ \pm 1 \} \times \sum^3$. The inverse images are points
or a nicely embedded $\sum^3$. Hence by Theorem 3.1,
$p_1$ is  a $UV^1$ map.

Next, the map 
$p_2 : M \vee M' \to M$ is the projection.
Since $\pi_1 (M') = \{ 1 \}$, Theorem 3.1 again implies
that $p_2$ is a $UV^1$ map.

It remains to observes that the composition of $UV^1$--maps
is again a $UV^1$--map.
\qed\enddemo

The proof of Lemma 3.2 shows the following:

\proclaim{Lemma 3.3}
If $M_1 \# M_2$ is a connected sum of
topological manifolds, then the smash map $p : M_1 \# M_2 \to M_1 \vee M_2$,
is $UV^1$.
\endproclaim

We now consider the following composition
$$
p = p_3 \circ p_2 \circ p_1 : 
(\mathop\#\limits _1\limits^r S^1 \times S^3) \# M' \;\;\; \mathop\to\limits^{p_1} \;\;\;
(\mathop\vee\limits _1\limits^r S^1 \times S^3) \vee M' \;\;\; 
$$
$$
\mathop\to\limits^{p_2} \;\;\;
(\mathop\vee\limits _1\limits^r S^1 \times S^3) \;\;\; \mathop\to\limits^{p_3} \;\;\;
\mathop\vee\limits _1\limits^r S^1 = B
$$
where $p_3$ is induced by the projection $S^1 \times S^3 \to S^1$.
Obviously, $p_3$ is a $UV^1$ map, hence
we obtain the following:

\proclaim{Corollary 3.4}
The map
$p : (\mathop\#\limits _1\limits^r S^1 \times S^3) \# M' = M \to
\mathop\vee\limits _1\limits^r S^1 = B$ is $UV^1$.
\endproclaim

From the Atiyah--Hirzebruch spectral sequence
$$
E_{rs}^2 = H_r (B, \pi_s (\Bbb L)) \Longrightarrow H_{r+s} (B, \Bbb L)
$$
we deduce the well--known fact
that $A : H_4 (B, \Bbb L) \to L_4 (\pi_1 (B))$ is an
isomorphism. In particular $E_{rs}^2 = 0$ for
$r > 1$, so the spectral sequence collapses and
$H_4 (B, \Bbb L) = H_0 (B, \Bbb Z) = \Bbb Z \cong L_4 (\pi_1(B))$.
Recall that:
$$
\pi_s (\Bbb L) =
\cases
0 &\text{if $s$ is odd} \\
\Bbb Z_2 &\text{if $s \equiv 2 (4)$} \\
\Bbb Z &\text{if $s \equiv 0 (4)$},
\endcases
$$
Since $M$ is a manifold, the map $p : M \to \mathop\vee\limits _1\limits^r S^1$ satisfies
conditions (i), (ii) and (iii). This proves Theorem 1.1. \qed

To prove Theorem 1.2 we consider
a Poincar\' e complex $X^4$ with $w_2 (X) = 0$ and
$\Lambda$--intersection form extended from $\Bbb Z$. The
fundamental group of $X$ is that of some surface $F$.
The construction of \cite{Cav-He-Rep} applies to give
a degree one normal map $X \mathop\to\limits^f F \times S^2$. This
splits the $\Lambda$--intersection form. Since it
is extended from the $\Bbb Z$--intersection form
one gets a homotopy equivalence $X \simeq F \times S^2 \# M' = M$,
where $M'$ is simply connected.

We get as above the following $UV^1$--map:
$$
p = p_3 \circ p_2 \circ p_1 : 
M \to F \times S^2 \vee M' \to F \times S^2 \to F = B.
$$

The Mayer--Vietoris technique can be applied to the
$L$--functor (see \cite{Capp}) and to $H_n (B, \Bbb L)$
to show that $A : H_4 (B, \Bbb L) \to L_4 (\pi_1(B))$
is an isomorphism. In particular
$L_4 (\pi_1 (B)) = \Bbb Z \oplus \Bbb Z_2$.
So $p : M \to B$ satisfies the conditions (i), (ii) and (iii),
which proves Theorem 1.2. \qed

Our next examples
are 4-manifolds whose fundamental groups
are knot groups. As the control space
$B$ we take a spine of the knot complement
in $S^3$. It is well known that $B$ is
an aspherical space and that $H_p (B) = \Bbb Z$
for $p = 0, 1$ and trivial otherwise. The Atiyah--Hirzebruch spectral sequence
then gives $H_4 (B, \Bbb L) = \Bbb Z$, in fact
$
A : H_4 (B, \Bbb L) \to L_4 (\pi_1 (B))
$
is an isomorphism.

\proclaim{Theorem 3.5}
Let $X = \partial (S^3 \setminus \mathop{N} \limits^o (k)) \times D^2$ be the boundary of
a regular neighborhood of the spine of the complement of
a torus knot, embedded in $\Bbb R^5$.
Then the surgery sequence
$$
\Cal S (X) \to [X, G/_{TOP}] \to L_4 (\pi_1 (X))
$$
is exact.
\endproclaim

Note that $X$ is a manifold. So it remains to verify
only condition (i) of $\S 1$. We state it as follows:

\proclaim{Lemma 3.6}
Let $p : X \to B$ be the restriction of
the neighborhood collapsing map. Then
$p$ is a $UV^1$ map.
\endproclaim

\demo{Proof}
We shall show
that the inverse images of points are 
$UV^1$  subsets of $X$ 
and then we shall apply Theorem 3.1 to get the assertion.
So let $k \subset S^1 \times S^1 \subset S^3$
be a torus knot in $S^3$ of type $(a, b)$, where $(a,b)=1$, and the
torus 
$S^1 \times S^1$
divides $S^3$ into two solid tori $T$ and $T^*$ such that
$k \subset T \cap T^* = \partial T = \partial T^* = S^1 \times S^1$.

Let $M^3 = S^3 \setminus \mathop{N^3} \limits^o (k)$,
where $N^{3} (k)$ is a small tubular neighborhood of the knot $k$ in $S^3$.
The spine of $M^3$ (i.e. a compact 2--polyhedron onto which the 3--manifold $M^3$ collapses) 
consists of $2$--dimensional compact polyhedra
$\Sigma \subset T$ and $\Sigma^* \subset T^*$, intersecting in 
an annulus $A=l \times [-1,1]$ which lies on
$T \cap T^* = S^1 \times S^1$ and where the curve $l$ is parallel on
$S^1 \times S^1$ to the knot $k$. (see Figure 1).

\centerline{
%
\beginpicture 
\setcoordinatesystem units <1.000mm,1.000mm>
\setplotarea x from -50 to 50, y from -30 to 30
\put {Figure 1} at 0 -25
\put {$A=l \times [-1,1]$} at 0 25
\setsolid
\setplotsymbol ({\fiverm.})
\plot -34.157 20.000 34.779 20.000 /
\plot -34.157 -20.000 34.779 -20.000 /
\setdashes
\plot -34.157 20.000 -50.000 20.000 /
\plot -34.157 -20.000 -50.000 -20.000 /
\plot 50.000 20.000 34.779 20.000 /
\plot 50.000 -20.000 34.779 -20.000 /
\setsolid
\ellipticalarc axes ratio 1:3 180 degrees from -34.157 20.000 center at -34.157 0.000
\ellipticalarc axes ratio 1:3 360 degrees from 34.157 20.000 center at 34.157 0.000
\circulararc 78.463 degrees from -5.112 17.642 center at -8.000 17.052
\circulararc 78.463 degrees from 10.888 17.642 center at 8.000 17.052
\circulararc 42.488 degrees from -23.202 -20.000 center at -23.243 -12.491
\circulararc 43.638 degrees from -7.202 -20.000 center at -7.202 -12.761
\plot 26.228 -6.000 27.700 -6.000 /
\plot 24.071 -8.000 27.900 -8.000 /
\plot 21.846 -10.000 28.000 -10.000 /
\plot 19.599 -12.000 28.700 -12.000 /
\plot 17.447 -14.000 29.300 -14.000 /
\plot 15.354 -16.000 30.000 -16.000 /
\plot 13.266 -18.000 28.266 -18.000 /
\plot -18.141 -18.000 -2.141 -18.000 /
\plot -15.818 -16.000 0.182 -16.000 /
\plot -13.844 -14.000 2.156 -14.000 /
\plot -12.147 -12.000 3.853 -12.000 /
\plot -10.679 -10.000 5.321 -10.000 /
\plot -9.408 -8.000 6.592 -8.000 /
\plot -8.311 -6.000 7.689 -6.000 /
\plot -7.370 -4.000 8.630 -4.000 /
\plot -6.572 -2.000 9.428 -2.000 /
\plot -5.908 0.000 10.092 0.000 /
\plot -5.370 2.000 10.630 2.000 /
\plot -4.951 4.000 11.049 4.000 /
\plot -4.649 6.000 11.351 6.000 /
\plot -4.460 8.000 11.540 8.000 /
\plot -4.383 10.000 11.617 10.000 /
\plot -4.416 12.000 11.584 12.000 /
\plot -4.560 14.000 11.440 14.000 /
\plot -4.816 16.000 11.184 16.000 /
\plot -5.209 18.000 10.813 18.000 /
\circulararc 64.667 degrees from -2.858 -18.552 center at -24.575 10.403
\circulararc 64.667 degrees from -18.858 -18.552 center at -40.575 10.403
\setplotsymbol ({\svtnrm.})
\circulararc 71.246 degrees from 2.791 18.000 center at -0.000 17.052
\circulararc 64.278 degrees from -10.206 -18.000 center at -32.359 10.344
\circulararc 41.870 degrees from -15.202 -20.000 center at -15.202 -12.761
\circulararc 41.765 degrees from 15.202 -20.000 center at 15.113 -12.363
\plot 20.266 -18.000 28.400 -11.0473 /
\setdashes
\setplotsymbol ({\fiverm.})
\ellipticalarc axes ratio 1:3 360 degrees from -34.157 -20.000 center at -34.157 0.000
\ellipticalarc axes ratio 1:3 360 degrees from 34.157 -20.000 center at 34.157 0.000
\circulararc 33.401 degrees from 19.217 -18.804 center at 23.202 -12.761
\circulararc 65.465 degrees from 5.112 17.642 center at 40.514 10.390
\circulararc 78.463 degrees from 8.000 20.000 center at 8.000 17.052
\circulararc 41.661 degrees from 2.858 -18.552 center at 7.332 -12.519
\circulararc 64.667 degrees from -10.888 17.642 center at 24.575 10.403
\circulararc 78.463 degrees from -8.000 20.000 center at -8.000 17.052
\setplotsymbol ({\svtnrm.})
\circulararc 36.870 degrees from 10.858 -18.552 center at 15.202 -12.761
\circulararc 64.667 degrees from -2.888 17.642 center at 32.575 10.403
\circulararc 78.463 degrees from 0.000 20.000 center at 0.000 17.052
\circulararc 41.765 degrees from -20.266 -18.000 center at -15.113 -12.363
\plot -20.266 -18.000 -28.000 -11.0473 /
\setsolid
\setplotsymbol ({\fiverm.})
\circulararc 41.765 degrees from -13.066 -18.000 center at -7.913 -12.363
\plot 13.066 -18.000 27.500 -4.805 /
\circulararc 41.765 degrees from 8.002 -20.000 center at 7.913 -12.363
\circulararc 41.765 degrees from -15.202 -20.000 center at -15.291 -12.363
\circulararc 41.765 degrees from 23.256 -20.000 center at 23.167 -12.363
\plot 28.320 -18.000 30.300 -16.190 /
\setdashes
\plot -13.066 -18.000 -27.500 -4.805 /
\plot -28.320 -18.000 -30.300 -16.190 /
\circulararc 41.765 degrees from -28.320 -18.000 center at -23.167 -12.363
\setshadesymbol ({\fiverm.})
\setshadegrid span <0.03in>
\vshade -40.823 0.000 0.000 <z,z,,>
-40 -9.630 9.630 <z,z,,>
-39 -13.744 13.744 <z,z,,>
-38 -16.343 16.343 <z,z,,>
-37 -18.090 18.090 <z,z,,>
-36 -19.220 19.220 <z,z,,>
-34.157 -20.000 20.000 <z,z,,>
27.492 -20.000 20.000 
28 -20.000 -7.670 <z,z,,>
29 -20.000 -12.675 <z,z,,>
30 -20.000 -15.646 <z,z,,>
31 -20.000 -17.615 <z,z,,>
32 -20.000 -18.924 <z,z,,>
34.157 -20.000 -20.000 /
\vshade 28 7.670 20.000 <z,z,,>
29 12.675 20.000 <z,z,,>
30 15.636 20.000 <z,z,,>
31 17.615 20.000 <z,z,,>
32 18.934 20.000 <z,z,,>
34.157 20.000 20.000 /
\setshadegrid span <0.02in>
\vshade 27.492 0.000 0.000 <z,z,,>
28 -7.670 7.670 <z,z,,>
29 -12.675 12.675 <z,z,,>
30 -15.646 15.636 <z,z,,>
31 -17.615 17.615 <z,z,,>
32 -18.924 18.934 <z,z,,>
34.157 -20.000 20.000 <z,z,,>
36 -19.220 19.220 <z,z,,>
37 -18.090 18.090 <z,z,,>
38 -16.343 16.343 <z,z,,>
39 -13.744 13.744 <z,z,,>
40 -9.630 9.630 <z,z,,>
40.823 0.000 0.000 /
\endpicture
%
}

If we look at the $2$--disk cross sections $C$ and $C^*$
of the solid tori $T$ and $T^*$ (which are
orthogonal to the longitudes of these solid tori),
the pictures of $\Sigma$ and $\Sigma^*$, respectively, are as in the Figure 2 below
(for the case when $(a,b) = (3,2)$, i.e. when $k$ is the $(3,2)$--torus knot):
where $\Sigma$ and $\Sigma^*$ are depicted by bold lines.

The collapsing map $\rho : M^3 \to B = \Sigma \cup \Sigma^*$ is described
on Figure 3:

\centerline{
%
\beginpicture 
\setcoordinatesystem units <1.300mm,1.300mm> 
\setplotarea x from -48 to 48, y from -25 to 25
\put{Figure 2} at 0 -20
%
\setcoordinatesystem units <0.600mm,0.600mm> point at -50 0
\setplotarea x from -50 to 50, y from -50 to 50
\setsolid
\setplotsymbol ({\fiverm.})
\circulararc 165.638 degrees from -9.922 38.750 center at 0.000 40.000
\plot 1.069 30.057 7.570 0.000 /
\plot 0.000 30.000 0.000 0.000 /
\plot 2.162 30.236 15.882 0.000 /
\plot 3.299 30.559 26.003 0.000 /
\plot 4.495 31.067 40.000 -0.000 /
\plot 5.748 31.816 37.353 14.310 /
\plot 7.024 32.882 30.586 25.778 /
\plot 8.245 34.340 22.286 33.217 /
\arrow <2mm> [.2,.4] from 0.000 30.000 to 0.000 15.000
\arrow <2mm> [.2,.4] from 1.069 30.057 to 4.228 15.452
\arrow <2mm> [.2,.4] from 2.162 30.236 to 8.264 16.787
\arrow <2mm> [.2,.4] from 3.299 30.559 to 11.927 18.946
\arrow <2mm> [.2,.4] from 4.495 31.067 to 15.052 21.830
\arrow <2mm> [.2,.4] from 5.748 31.816 to 17.495 25.309
\arrow <2mm> [.2,.4] from 7.024 32.882 to 19.149 29.226
\arrow <2mm> [.2,.4] from 8.245 34.340 to 14.952 33.804
\arrow <2mm> [.2,.4] from -8.245 34.340 to -14.952 33.804
\arrow <2mm> [.2,.4] from -7.024 32.882 to -19.149 29.226
\arrow <2mm> [.2,.4] from -5.748 31.816 to -17.495 25.309
\arrow <2mm> [.2,.4] from -4.495 31.067 to -15.052 21.830
\arrow <2mm> [.2,.4] from -3.299 30.559 to -11.927 18.946
\arrow <2mm> [.2,.4] from -2.162 30.236 to -8.264 16.787
\arrow <2mm> [.2,.4] from -1.069 30.057 to -4.228 15.452
\plot -8.245 34.340 -22.286 33.217 /
\plot -7.024 32.882 -30.586 25.778 /
\plot -5.748 31.816 -37.353 14.310 /
\plot -4.495 31.067 -40.000 0.000 /
\plot -3.299 30.559 -26.003 0.000 /
\plot -2.162 30.236 -15.882 -0.000 /
\plot -1.069 30.057 -7.570 -0.000 /
\plot -1.069 -30.057 -7.570 -0.000 /
\plot -2.162 -30.236 -15.882 -0.000 /
\plot -3.299 -30.559 -26.003 0.000 /
\plot -4.495 -31.067 -40.000 0.000 /
\plot -5.748 -31.816 -37.353 -14.310 /
\plot -7.024 -32.882 -30.586 -25.778 /
\plot -8.245 -34.340 -22.286 -33.217 /
\arrow <2mm> [.2,.4] from -1.069 -30.057 to -4.228 -15.452
\arrow <2mm> [.2,.4] from -2.162 -30.236 to -8.264 -16.787
\arrow <2mm> [.2,.4] from -3.299 -30.559 to -11.927 -18.946
\arrow <2mm> [.2,.4] from -4.495 -31.067 to -15.052 -21.830
\arrow <2mm> [.2,.4] from -5.748 -31.816 to -17.495 -25.309
\arrow <2mm> [.2,.4] from -7.024 -32.882 to -19.149 -29.226
\arrow <2mm> [.2,.4] from -8.245 -34.340 to -14.952 -33.804
\arrow <2mm> [.2,.4] from 8.245 -34.340 to 14.952 -33.804
\arrow <2mm> [.2,.4] from 7.024 -32.882 to 19.149 -29.226
\arrow <2mm> [.2,.4] from 5.748 -31.816 to 17.495 -25.309
\arrow <2mm> [.2,.4] from 4.495 -31.067 to 15.052 -21.830
\arrow <2mm> [.2,.4] from 3.299 -30.559 to 11.927 -18.946
\arrow <2mm> [.2,.4] from 2.162 -30.236 to 8.264 -16.787
\arrow <2mm> [.2,.4] from 1.069 -30.057 to 4.228 -15.452
\arrow <2mm> [.2,.4] from 0.000 -30.000 to 0.000 -15.000
\plot 8.245 -34.340 22.286 -33.217 /
\plot 7.024 -32.882 30.586 -25.778 /
\plot 5.748 -31.816 37.353 -14.310 /
\plot 4.495 -31.067 40.000 0.000 /
\plot 3.299 -30.559 26.003 0.000 /
\plot 2.162 -30.236 15.882 0.000 /
\plot 0.000 -30.000 0.000 0.000 /
\plot 1.069 -30.057 7.570 0.000 /
\circulararc 165.638 degrees from 9.922 -38.750 center at -0.000 -40.000
\setplotsymbol ({\svtnrm.})
\plot -40.000 0.000 40.000 0.000 /
\circulararc 151.276 degrees from 9.922 -38.750 center at 0.000 0.000
\circulararc 151.276 degrees from -9.922 38.750 center at -0.000 -0.000
\put {$k$} at 0.000 -48.000
\put {$k$} at 0.000 48.000
\put {$\bullet$} at 0.000 40.000
\circulararc 360 degrees from 1.000 40.000 center at 0.000 40.000
\circulararc 360 degrees from 2.000 40.000 center at 0.000 40.000
\circulararc 360 degrees from 3.000 40.000 center at 0.000 40.000
\put {$\bullet$} at -0.000 -40.000
\circulararc 360 degrees from 1.000 -40.000 center at 0.000 -40.000
\circulararc 360 degrees from 2.000 -40.000 center at 0.000 -40.000
\circulararc 360 degrees from 3.000 -40.000 center at 0.000 -40.000
\put {$a^*$} at 3.000 -5.000
\put {$\bullet$} at 0.000 -0.000
\put {$b^*$} at 28.550 -5.000
\put {$\bullet$} at 25.550 0.000
\put {$C^*$} at -35.000 35.000
\put {$c^*$} at 45.368 15.096
\put {$d^*$} at 15.078 42.719
\put {$\bullet$} at 9.922 38.750
\put {$\bullet$} at 37.353 14.310
\setcoordinatesystem units <0.600mm,0.600mm> point at 50 0
\setplotarea x from -50 to 50, y from -50 to 50
\setsolid
\setplotsymbol ({\fiverm.})
\circulararc 165.638 degrees from -9.922 38.750 center at -0.000 40.000
\plot 1.463 30.108 8.924 5.153 /
\plot 2.984 30.456 16.666 9.622 /
\plot 4.608 31.125 24.679 14.248 /
\plot 6.330 32.259 34.641 20.000 /
\plot 8.035 34.047 23.756 32.181 /
\arrow <2mm> [.2,.4] from 8.035 34.047 to 18.654 32.787
\arrow <2mm> [.2,.4] from 6.330 32.259 to 20.729 26.024
\arrow <2mm> [.2,.4] from 4.608 31.125 to 16.446 21.171
\arrow <2mm> [.2,.4] from 1.463 30.108 to 5.778 15.677
\arrow <2mm> [.2,.4] from 2.984 30.456 to 11.346 17.723
\arrow <2mm> [.2,.4] from -2.984 30.456 to -11.346 17.723
\arrow <2mm> [.2,.4] from -1.463 30.108 to -5.778 15.677
\arrow <2mm> [.2,.4] from -0.000 30.000 to 0.000 15.000
\arrow <2mm> [.2,.4] from -4.608 31.125 to -16.446 21.171
\arrow <2mm> [.2,.4] from -6.330 32.259 to -20.729 26.024
\arrow <2mm> [.2,.4] from -8.035 34.047 to -18.654 32.787
\plot -8.035 34.047 -23.756 32.181 /
\plot -6.330 32.259 -34.641 20.000 /
\plot -4.608 31.125 -24.679 14.248 /
\plot -2.984 30.456 -16.666 9.622 /
\plot -1.463 30.108 -8.924 5.153 /
\plot -0.000 30.000 -0.000 0.000 /
\plot -0.000 0.000 -0.000 -40.000 /
\plot -25.981 -15.000 -0.000 0.000 /
\plot -26.805 -13.787 -8.924 5.153 /
\plot -27.868 -12.643 -16.666 9.622 /
\plot -29.259 -11.572 -24.679 14.248 /
\plot -31.102 -10.647 -34.641 20.000 /
\plot -33.503 -10.065 -39.748 4.483 /
\arrow <2mm> [.2,.4] from -33.503 -10.065 to -37.721 -0.239
\arrow <2mm> [.2,.4] from -31.102 -10.647 to -32.902 4.939
\arrow <2mm> [.2,.4] from -29.259 -11.572 to -26.558 3.657
\arrow <2mm> [.2,.4] from -25.981 -15.000 to -12.990 -7.500
\arrow <2mm> [.2,.4] from -26.805 -13.787 to -16.465 -2.835
\arrow <2mm> [.2,.4] from -27.868 -12.643 to -21.022 0.964
\arrow <2mm> [.2,.4] from -24.883 -17.812 to -9.675 -18.687
\arrow <2mm> [.2,.4] from -25.343 -16.321 to -10.688 -12.842
\arrow <2mm> [.2,.4] from -24.651 -19.553 to -10.112 -24.828
\arrow <2mm> [.2,.4] from -24.772 -21.612 to -12.173 -30.964
\arrow <2mm> [.2,.4] from -25.468 -23.982 to -19.067 -32.548
\plot -25.468 -23.982 -15.992 -36.664 /
\plot -24.772 -21.612 0.000 -40.000 /
\plot -24.651 -19.553 -0.000 -28.496 /
\plot -24.883 -17.812 -0.000 -19.244 /
\plot -25.343 -16.321 -0.000 -10.305 /
\circulararc 165.638 degrees from -28.598 -27.967 center at -34.641 -20.000
\circulararc 165.638 degrees from 38.519 -10.783 center at 34.641 -20.000
\plot 25.343 -16.321 0.000 -10.305 /
\plot 24.883 -17.812 0.000 -19.244 /
\plot 24.651 -19.553 0.000 -28.496 /
\plot 24.772 -21.612 0.000 -40.000 /
\plot 25.468 -23.982 15.992 -36.664 /
\arrow <2mm> [.2,.4] from 25.468 -23.982 to 19.067 -32.548
\arrow <2mm> [.2,.4] from 24.772 -21.612 to 12.173 -30.964
\arrow <2mm> [.2,.4] from 24.651 -19.553 to 10.112 -24.828
\arrow <2mm> [.2,.4] from 25.343 -16.321 to 10.688 -12.842
\arrow <2mm> [.2,.4] from 24.883 -17.812 to 9.675 -18.687
\arrow <2mm> [.2,.4] from 27.868 -12.643 to 21.022 0.964
\arrow <2mm> [.2,.4] from 26.805 -13.787 to 16.465 -2.835
\arrow <2mm> [.2,.4] from 25.981 -15.000 to 12.990 -7.500
\arrow <2mm> [.2,.4] from 29.259 -11.572 to 26.558 3.657
\arrow <2mm> [.2,.4] from 31.102 -10.647 to 32.902 4.939
\arrow <2mm> [.2,.4] from 33.503 -10.065 to 37.721 -0.239
\plot 33.503 -10.065 39.748 4.483 /
\plot 31.102 -10.647 34.641 20.000 /
\plot 29.259 -11.572 24.679 14.248 /
\plot 27.868 -12.643 16.666 9.622 /
\plot 26.805 -13.787 8.924 5.153 /
\plot 25.981 -15.000 0.000 0.000 /
\setplotsymbol ({\svtnrm.})
\circulararc 91.277 degrees from -28.598 -27.967 center at 0.000 -0.000
\circulararc 91.277 degrees from -9.922 38.750 center at 0.000 -0.000
\circulararc 91.277 degrees from 38.519 -10.783 center at -0.000 -0.000
\plot 0.000 0.000 0.000 -40.000 /
\plot 0.000 0.000 -34.641 20.000 /
\plot 34.641 20.000 0.000 0.000 /
\put {$b$} at 30.278 13.383
\put {$a$} at 6.990 -0.321
\put {$c$} at 29.061 33.726
\put {$k$} at 39.239 -26.500
\put {$k$} at 0.000 48.000
\put {$k$} at -39.239 -26.500
\put {$C$} at -35.000 35.000
\put {$d$} at 15.078 42.719
\put {$\bullet$} at 23.757 32.181
\put {$\bullet$} at 9.922 38.750
\put {$\bullet$} at 24.679 14.248
\put {$\bullet$} at -0.000 -0.000
\put {$\bullet$} at -0.000 40.000
\circulararc 360 degrees from 1.000 40.000 center at 0.000 40.000
\circulararc 360 degrees from 2.000 40.000 center at 0.000 40.000
\circulararc 360 degrees from 3.000 40.000 center at 0.000 40.000
\put {$\bullet$} at -34.641 -20.000
\circulararc 360 degrees from -35.641 -20.000 center at -34.641 -20.000
\circulararc 360 degrees from -36.641 -20.000 center at -34.641 -20.000
\circulararc 360 degrees from -37.641 -20.000 center at -34.641 -20.000
\put {$\bullet$} at 34.641 -20.000
\circulararc 360 degrees from 35.641 -20.000 center at 34.641 -20.000
\circulararc 360 degrees from 36.641 -20.000 center at 34.641 -20.000
\circulararc 360 degrees from 37.641 -20.000 center at 34.641 -20.000
\endpicture
%
}

\centerline{
%
\beginpicture 
\setcoordinatesystem units <1.300mm,1.300mm> 
\setplotarea x from -48 to 48, y from -25 to 25
\put{Figure 3} at 0 -20
%
\setcoordinatesystem units <0.600mm,0.600mm> point at -50 0
\setplotarea x from -50 to 50, y from -50 to 50
\setsolid
\setplotsymbol ({\svtnrm.})
\plot -40.000 0.000 40.000 0.000 /
\circulararc 151.276 degrees from 9.922 -38.750 center at 0.000 0.000
\circulararc 151.276 degrees from -9.922 38.750 center at -0.000 -0.000
\put {$k$} at 0.000 -48.000
\put {$k$} at 0.000 48.000
\put {$\bullet$} at 0.000 40.000
\circulararc 360 degrees from 1.000 40.000 center at 0.000 40.000
\circulararc 360 degrees from 2.000 40.000 center at 0.000 40.000
\circulararc 360 degrees from 3.000 40.000 center at 0.000 40.000
\put {$\bullet$} at -0.000 -40.000
\circulararc 360 degrees from 1.000 -40.000 center at 0.000 -40.000
\circulararc 360 degrees from 2.000 -40.000 center at 0.000 -40.000
\circulararc 360 degrees from 3.000 -40.000 center at 0.000 -40.000
\put {$a^*$} at 3.000 -5.000
\put {$\bullet$} at 0.000 -0.000
\put {$b^*$} at 28.550 -5.000
\put {$\bullet$} at 25.550 0.000
\put {$\Sigma^*$} at -35.000 35.000
\put {$c^*$} at 45.368 15.096
\put {$d^*$} at 15.078 42.719
\put {$\bullet$} at 9.922 38.750
\put {$\bullet$} at 37.353 14.310
\setcoordinatesystem units <0.600mm,0.600mm> point at 50 0
\setplotarea x from -50 to 50, y from -50 to 50
\setsolid
\setplotsymbol ({\fiverm.})
\setplotsymbol ({\svtnrm.})
\circulararc 91.277 degrees from -28.598 -27.967 center at 0.000 -0.000
\circulararc 91.277 degrees from -9.922 38.750 center at 0.000 -0.000
\circulararc 91.277 degrees from 38.519 -10.783 center at -0.000 -0.000
\plot 0.000 0.000 0.000 -40.000 /
\plot 0.000 0.000 -34.641 20.000 /
\plot 34.641 20.000 0.000 0.000 /
\put {$b$} at 30.278 13.383
\put {$a$} at 6.990 -0.321
\put {$c$} at 29.061 33.726
\put {$k$} at 39.239 -26.500
\put {$k$} at 0.000 48.000
\put {$k$} at -39.239 -26.500
\put {$\Sigma$} at -35.000 35.000
\put {$d$} at 15.078 42.719
\put {$\bullet$} at 23.757 32.181
\put {$\bullet$} at 9.922 38.750
\put {$\bullet$} at 24.679 14.248
\put {$\bullet$} at -0.000 -0.000
\put {$\bullet$} at -0.000 40.000
\circulararc 360 degrees from 1.000 40.000 center at 0.000 40.000
\circulararc 360 degrees from 2.000 40.000 center at 0.000 40.000
\circulararc 360 degrees from 3.000 40.000 center at 0.000 40.000
\put {$\bullet$} at -34.641 -20.000
\circulararc 360 degrees from -35.641 -20.000 center at -34.641 -20.000
\circulararc 360 degrees from -36.641 -20.000 center at -34.641 -20.000
\circulararc 360 degrees from -37.641 -20.000 center at -34.641 -20.000
\put {$\bullet$} at 34.641 -20.000
\circulararc 360 degrees from 35.641 -20.000 center at 34.641 -20.000
\circulararc 360 degrees from 36.641 -20.000 center at 34.641 -20.000
\circulararc 360 degrees from 37.641 -20.000 center at 34.641 -20.000
\endpicture
%
}
Let us consider the point inverses of the map $\rho$. There are essentially
three different types of points to consider:
$\rho^{-1} (a)$ (resp. $\rho^{-1} (a^*)$)
is a bouquet of 3 (resp. 2) intervals, 
$\rho^{-1} (b)$ and $\rho^{-1} (c)$ (resp. $\rho^{-1} (b^*)$ and $\rho^{-1} (c^*)$)
are bouquets of 2 intervals, and finally,
$\rho^{-1} (d)$ (resp. $\rho^{-1} (d^*)$) is just a point.

We now consider the embedding $B \subset \Bbb R^5$ given
by 
$$B \subset M^3 \subset S^3 \subset S^3 \times D^1 \subset \Bbb R^4 \times \{0\} \subset \Bbb R^5.$$
Let $X$ be the boundary of the regular neighborhood $M^3 \times D^2$
of
$B \subset \Bbb R^5$:
$$X = \partial (M^3 \times D^2) = \partial M^3 \times D^2 \cup M^3 \times \partial D^2.$$
Then the collapsing map 
$p  : X \to B$
is
the composition of the canonical projection 
$\pi: M^3 \times D^2 \to M^3$
followed by the collapsing map 
$\rho : M^3 \to B$ described above,
that is
$p = \rho \circ (\pi |_X).$

It can now be easily verified that $p^{-1} (z) = (\pi |_X) ^{-1}(\rho^{-1}(z))$ is indeed 
a $UV^1$ subset of $X$, for every point $z \in B$.
We shall do this for interior points $a$ and $b$ (resp. $a^*$ and $b^*$)
of $\Sigma$ (resp.  $\Sigma^*$) and we leave the reader the verification for the boundary points
$c$ and $d$ (resp. $c^*$ and $d^*$).
As we have seen above, in such a case, $\rho^{-1}(z)$ is a bouquet of arcs 
$\mathop\vee\limits _{i=1} \limits^m D_i^1 $ with endpoints $d_{i}$. 
Therefore we can easily see that 
$p^{-1} (z)$ is in this case a 2-sphere 
with finitely many disks
attached along the equator (see Figure 4).
Therefore 
$p^{-1} (z)$ is certainly a $UV^1$ subset:
$$
\Bbb R^3 \times \Bbb R^2 \supset \partial (M^3 \times D^2) \supset 
p^{-1} (z) =
\partial (\mathop\vee\limits _{i=1} \limits^m D_i^1 \times D^2) =
$$
$$
= (\partial (\mathop\vee\limits _{i=1} \limits^m D_i^1) \times D^2) \cup 
(\mathop\vee\limits _{i=1} \limits^m D_i^1 \times \partial D^2) 
\simeq \mathop\vee\limits _{i=1} \limits^{m-1} S^2.
$$

\centerline{
%
\beginpicture
\setcoordinatesystem units <1.000mm,1.000mm>
\setplotarea x from -50 to 50, y from -25 to 25
\put {Figure 4} at 0.000 -23.000
\setcoordinatesystem units <1.000mm,1.000mm> point at 25 0
\setplotarea x from -25 to 25, y from -25 to 25
\setsolid
\setplotsymbol ({\fiverm.})
\circulararc 360.000 degrees from 20.000 0.000 center at 0.000 0.000
\ellipticalarc axes ratio 5:1 180.000 degrees from -20.000 0.000 center at 0.000 0.000
\ellipticalarc axes ratio 5:1 180.000 degrees from -10.500 16.500 center at 0.000 16.500
\ellipticalarc axes ratio 5:1 180.000 degrees from -10.500 -16.500 center at 0.000 -16.500
\ellipticalarc axes ratio 5:1 180.000 degrees from -10.500 0.000 center at 0.000 0.000
\ellipticalarc axes ratio 5:1 180.000 degrees from 10.500 0.000 center at 0.000 0.000
\setdashes
\ellipticalarc axes ratio 5:1 180.000 degrees from 20.000 0.000 center at 0.000 0.000
\ellipticalarc axes ratio 5:1 180.000 degrees from 10.500 16.500 center at 0.000 16.500
\ellipticalarc axes ratio 5:1 180.000 degrees from 10.500 -16.500 center at 0.000 -16.500
\put {$D_1$} at 0.000 -23.000
\put {$D_2$} at 5.000 0.000
\put {$D_3$} at 0.000 23.000
\put {$A_1$} at -13.000 -10.000
\put {$A_2$} at -13.000 0.000
\put {$A_3$} at -13.000 10.000
\put {$Q$} at 18.000 18.000
\setshadesymbol ({\fiverm.})
\setshadegrid span <0.02in>
\vshade -10.500 0.000 0.000 <z,z,,>
-10 -0.640 0.640 <z,z,,>
-9 -1.081 1.081 <z,z,,>
-7 -1.565 1.565 <z,z,,>
-5 -1.847 1.847 <z,z,,>
-3 -2.012 2.012 <z,z,,>
0 -2.100 2.100 <z,z,,>
3 -2.012 2.012 <z,z,,>
5 -1.847 1.847 <z,z,,>
7 -1.565 1.565 <z,z,,>
9 -1.081 1.081 <z,z,,>
10 -0.640 0.640 <z,z,,>
10.500 0.000 0.000 /
\vshade -10.500 16.500 16.500 <z,z,,>
-10 15.859 17.320 <z,z,,>
-9 15.418 17.861 <z,z,,>
-7 14.934 18.735 <z,z,,>
-5 14.653 19.365 <z,z,,>
-3 14.487 19.773 <z,z,,>
0 14.400 20.000 <z,z,,>
3 14.487 19.773 <z,z,,>
5 14.653 19.365 <z,z,,>
7 14.934 18.735 <z,z,,>
9 15.418 17.861 <z,z,,>
10 15.859 17.320 <z,z,,>
10.500 16.500 16.500 /
\vshade -10.500 -16.500 -16.500 <z,z,,>
-10 -17.320 -17.140 <z,z,,>
-9 -17.861 -17.581 <z,z,,>
-7 -18.735 -18.065 <z,z,,>
-5 -19.365 -18.347 <z,z,,>
-3 -19.773 -18.512 <z,z,,>
0 -20.000 -18.600 <z,z,,>
3 -19.773 -18.512 <z,z,,>
5 -19.365 -18.347 <z,z,,>
7 -18.735 -18.065 <z,z,,>
9 -17.861 -17.581 <z,z,,>
10 -17.320 -17.140 <z,z,,>
10.500 -16.500 -16.500 /
\setshadegrid span <0.015in>
\vshade -20.000 0.000 0.000 <z,z,,>
-19 -1.249 1.249 <z,z,,>
-17 -2.107 2.107 <z,z,,>
-16 -2.400 2.400 <z,z,,>
-14 -2.857 2.857 <z,z,,>
-12 -3.200 3.200 <z,z,,>
-10.500 -3.404 3.404 /
\vshade -10.500 0.000 3.404 <z,z,,>
-10 0.640 3.464 <z,z,,>
-9 1.081 3.572 <z,z,,>
-7 1.565 3.747 <z,z,,>
-5 1.847 3.873 <z,z,,>
-3 2.012 3.955 <z,z,,>
0 2.100 4.000 <z,z,,>
3 2.012 3.955 <z,z,,>
5 1.847 3.873 <z,z,,>
7 1.565 3.747 <z,z,,>
9 1.081 3.572 <z,z,,>
10 0.640 3.464 <z,z,,>
10.500 0.000 3.404 /
\vshade -10.500 -3.404 0.000 <z,z,,>
-10 -3.464 -0.640 <z,z,,>
-9 -3.572 -1.081 <z,z,,>
-7 -3.747 -1.565 <z,z,,>
-5 -3.873 -1.847 <z,z,,>
-3 -3.955 -2.012 <z,z,,>
0 -4.000 -2.100 <z,z,,>
3 -3.955 -2.012 <z,z,,>
5 -3.873 -1.847 <z,z,,>
7 -3.747 -1.565 <z,z,,>
9 -3.572 -1.081 <z,z,,>
10 -3.464 -0.640 <z,z,,>
10.500 -3.404 0.000 /
\vshade 10.500 -3.404 3.404 <z,z,,>
12 -3.200 3.200 <z,z,,>
14 -2.857 2.857 <z,z,,>
16 -2.400 2.400 <z,z,,>
17 -2.107 2.107 <z,z,,>
19 -1.249 1.249 <z,z,,>
20.000 0.000 0.000 /
\setcoordinatesystem units <1.000mm,1.000mm> point at -25 0
\setplotarea x from -25 to 25, y from -25 to 25
\setsolid
\setplotsymbol ({\fiverm.})
\circulararc 360.000 degrees from 20.000 0.000 center at 0.000 0.000
\ellipticalarc axes ratio 5:1 180.000 degrees from -20.000 0.000 center at 0.000 0.000
\ellipticalarc axes ratio 5:1 180.000 degrees from -10.500 16.500 center at 0.000 16.500
\ellipticalarc axes ratio 5:1 180.000 degrees from -10.500 -16.500 center at 0.000 -16.500
\setdashes
\ellipticalarc axes ratio 5:1 180.000 degrees from 20.000 0.000 center at 0.000 0.000
\ellipticalarc axes ratio 5:1 180.000 degrees from 10.500 16.500 center at 0.000 16.500
\ellipticalarc axes ratio 5:1 180.000 degrees from 10.500 -16.500 center at 0.000 -16.500
\put {$Q^*$} at 18.000 18.000
\put {$D_1^*$} at 0.000 -23.000
\put {$D_2^*$} at 0.000 23.000
\put {$A_1^*$} at -13.000 -10.000
\put {$A_2^*$} at -13.000 10.000
\setshadegrid span <0.02in>
\vshade -10.500 16.500 16.500 <z,z,,>
-10 15.859 17.320 <z,z,,>
-9 15.418 17.861 <z,z,,>
-7 14.934 18.735 <z,z,,>
-5 14.653 19.365 <z,z,,>
-3 14.487 19.773 <z,z,,>
0 14.400 20.000 <z,z,,>
3 14.487 19.773 <z,z,,>
5 14.653 19.365 <z,z,,>
7 14.934 18.735 <z,z,,>
9 15.418 17.861 <z,z,,>
10 15.859 17.320 <z,z,,>
10.500 16.500 16.500 /
\vshade -10.500 -16.500 -16.500 <z,z,,>
-10 -17.320 -17.140 <z,z,,>
-9 -17.861 -17.581 <z,z,,>
-7 -18.735 -18.065 <z,z,,>
-5 -19.365 -18.347 <z,z,,>
-3 -19.773 -18.512 <z,z,,>
0 -20.000 -18.600 <z,z,,>
3 -19.773 -18.512 <z,z,,>
5 -19.365 -18.347 <z,z,,>
7 -18.735 -18.065 <z,z,,>
9 -17.861 -17.581 <z,z,,>
10 -17.320 -17.140 <z,z,,>
10.500 -16.500 -16.500 /
\endpicture
%
}
It therefore follows by Theorem 3.2 that
$p$ is indeed
a $UV^1$--map, as it was asserted.
\qed\enddemo

{\sl Remark.} Yamasaki has recently proved 
that our strategy
can also be used to 
prove Theorem 3.5 for hyperbolic knots 
\cite{Ya}.
 
\head {$\S$ 4. Appendix} \endhead

\noindent
{\bf 1.}
Using Daverman's theorem (see Theorem 3.1 above)
one gets $UV^1$--maps $p : X \to B$ in the following
way: Let $U \subset X$ be a compact $4$--dimensional 
submanifold
with boundary $\partial U$. If every component of $U$ is simply connected
then the projection 
$p_U : X \to X/_{U}$ is $UV^1$.
We shall say that 
a $4$--manifold $X$ is a {\sl good manifold}
if there exists such a submanifold $U$ that
$A : H_4 (X/_{U}; \Bbb L) \to L_4 (\pi_1 (X))$
is injective. Note that $\pi_1 (X) \cong \pi_1 (X/_{U})$.
Then one gets the following more general result:

\proclaim{Theorem 4.1}
If a 4--manifold 
$M^4$ is homotopy equivalent to a good
$4$--manifold $X$, then the surgery sequence
$$
\Cal S (X) \to [X, G/_{TOP}] \to L_4 (\pi_1)
$$
is exact.
\endproclaim

Note that "$X$ has a good fundamental
group" does not imply that "$X$
is a good manifold", and vice--versa.

\bigskip

\noindent
{\bf 2.}
Let $\pi$ be the fundamental group of an arbitrary knot $k \subset S^3$.
There is a well-known procedure (see e.g. \cite{Ma}) by which one can
construct a (special) compact $2$--polyhedron $K \subset S^3$,
such that $\pi_1 (K) \cong \pi = \pi_1 (\overline{S^3 \setminus k})$.
We briefly outline it below.

Let a knot $k$ be given by its projection onto $S^2$,
which is discontinued at the double points -- in order
to show which of the diagram's parts goes over the other. Glue a long closed
strip (a tunnel) to $S^2$ and to itself along the knot projection, as it is shown in
Figure 5.

%
\beginpicture 
\setcoordinatesystem units <0.500mm,0.500mm>
\setplotarea x from -120 to 130, y from -70 to 60
\put{Figure 5} at 0 -60

\setcoordinatesystem units <0.500mm,0.500mm> point at 90 20
\setplotarea x from -30 to 30, y from -30 to 30
\setsolid
\setplotsymbol ({\fiverm.})
\plot -4.207 14.668 -8.708 8.030 /
\plot 3.425 12.313 -4.207 14.668 /
\plot -1.468 5.460 3.425 12.313 /
\plot -8.708 8.030 -1.468 5.460 /
\circulararc 360 degrees from 25.000 0.000 center at 0.000 0.000
\ellipticalarc axes ratio 4:1  -180 degrees from 25.000 0.000 center at 0.000 0.000
\setdashes
\ellipticalarc axes ratio 4:1  360 degrees from 25.000 0.000 center at 0.000 0.000
\setsolid
\setplotsymbol ({\rm.})
\circulararc 4.099 degrees from -12.247 -7.071 center at -5.176 19.319
\circulararc 4.099 degrees from 12.247 -7.071 center at -14.142 -14.142
\circulararc 4.099 degrees from 0.000 14.142 center at 19.319 -5.176
\circulararc 45.411 degrees from -4.243 8.653 center at 19.319 -5.176
\circulararc 45.411 degrees from -5.372 -8.001 center at -5.176 19.319
\circulararc 45.411 degrees from 9.616 -0.652 center at -14.142 -14.142
\circulararc 180.000 degrees from 14.142 -0.000 center at 7.071 7.071
\circulararc 180.000 degrees from -7.071 12.247 center at -9.659 2.588
\circulararc 180.000 degrees from -7.071 -12.247 center at 2.588 -9.659

\circulararc -90.000 degrees from 0 20 center at 30 20
\plot 30 50 70 50 /
\plot 70 50 67 53 /
\plot 70 50 67 47 /

\setcoordinatesystem units <0.500mm,0.500mm> point at -50 0
\setplotarea x from -90 to 80, y from -65 to 60


\setsolid
\setplotsymbol ({\fiverm.})
\plot -20.000 50.710 -40.905 29.805 /
\plot 19.393 34.393 -20.000 50.710 /
\circulararc 22.680 degrees from -52.284 -20.770 center at -84.755 -99.164
\plot -52.284 -20.770 -36.036 -27.500 /
\plot -55.269 15.441 -85.022 -14.311 /
\circulararc 19.287 degrees from 76.355 5.645 center at 19.812 -57.624
\plot 76.355 5.645 55.373 -15.337 /
\plot 52.284 20.770 36.036 27.500 /
\plot -7.716 -39.230 -23.964 -32.500 /
\plot 44.731 -25.980 41.213 -29.497 /
\circulararc 24.412 degrees from 21.780 -59.658 center at -40.188 -117.624
\circulararc 24.412 degrees from 41.213 -29.497 center at 101.213 -89.497
\circulararc 6.372 degrees from 22.174 -49.355 center at 22.271 -49.439
\plot -15.876 -0.876 -55.269 15.441 /
\plot 44.731 -25.980 0.766 -7.769 /
\plot 16.802 8.267 36.036 27.500 /
\plot -23.964 -32.500 0.766 -7.769 /
\plot 0.735 15.735 25.000 40.000 /
\plot -35.000 -20.000 -10.607 4.393 /
\plot -28.497 -19.961 -36.036 -27.500 /
\plot 46.618 -11.711 55.373 -15.337 /
\circulararc 45.000 degrees from -7.769 5.766 center at -7.769 2.231
\circulararc 29.992 degrees from -47.680 32.094 center at -49.118 27.839
\circulararc 29.992 degrees from -47.680 32.094 center at -49.118 27.839
\circulararc 45.000 degrees from 27.500 41.036 center at 27.500 37.500
\circulararc 7.625 degrees from 4.013 14.586 center at 2.024 10.559
\plot 3.461 14.814 -47.680 32.094 /
\plot 52.871 -9.555 4.013 14.586 /
\circulararc 67.500 degrees from -50.000 32.243 center at -46.734 24.357
\plot -55.269 24.357 -55.269 15.441 /
\plot 44.731 -25.980 55.373 -15.337 /
\plot 44.731 -17.063 44.731 -25.980 /
\plot 55.373 -13.582 55.373 -15.337 /
\circulararc 101.323 degrees from 55.373 -13.582 center at 50.882 -13.582
\circulararc 67.500 degrees from 50.000 -9.178 center at 53.266 -17.063
\circulararc 67.646 degrees from 54.619 -11.091 center at 50.882 -13.582
\circulararc 90.000 degrees from 0.766 -2.769 center at -7.769 -2.769
\plot 0.766 -2.769 0.766 -7.769 /
\plot -36.036 -27.500 -23.964 -32.500 /
\plot 36.036 32.500 36.036 27.500 /
\plot -36.036 -22.500 -36.036 -27.500 /
\plot -23.964 -27.500 -23.964 -32.500 /
\circulararc 90.000 degrees from 36.036 32.500 center at 27.500 32.500
\circulararc 90.000 degrees from -23.964 -27.500 center at -32.500 -27.500
\circulararc 90.000 degrees from -32.500 -18.964 center at -32.500 -22.500
\setdots
\setplotsymbol ({\rm.})
\plot 9.192 9.192 30.000 30.000 /
\plot -12.071 5.000 -50.000 20.710 /
\circulararc 90.000 degrees from -5.000 12.071 center at -5.000 5.000
\circulararc 90.000 degrees from 12.071 -5.000 center at -5.000 -5.000
\plot 50.000 -20.710 12.071 -5.000 /
\plot -30.000 -30.000 4.950 4.950 /
\endpicture
%

We obtain a (special) spine $K$ of the twice punctured knot complement, i.e.
of the compact $3$--dimensional manifold $N^3 \subset S^3$ which collapses to
$K$, $N^3 \searrow K$, and $\partial N = S^2 \cup (S^1 \times S^1) \cup S^2$.
Unfortunately, the $2$--polyhedron $K$ contains a nontrivial
$2$--sphere, i.e. $H_2 (K, \Bbb Z_2) \ne 0$.
However, since $L_4 (\pi) = \Bbb Z$, the control space
$B = K$ should have trivial $H_2 (K, \Bbb Z)$,
in order to guarantee that the assembly map
$A : H_4 (B, \Bbb L) \to L_4 (\pi)$ is an isomorphism.

\head {Acknowledgements} \endhead
The second author was supported in part by the
Ministry for Education, Science and Sport of the
Republic of Slovenia. This paper is a result of our {\sl Research in Pairs}
stay at the Mathematisches Forschungsinstitut Oberwolfach in April 2004. 
We thank the institute for the excellent working conditions and their hospitality. 
We also acknowledge the referee for several important
remarks and suggestions.

\vfill\eject

\enddocument
\end